\documentclass[conference]{IEEEtran}
\usepackage{mwe} 
\usepackage[nopar]{kantlipsum} 
\usepackage{amsmath}
\usepackage{amssymb}
\usepackage{color,soul}
\usepackage{graphicx}
\usepackage{epstopdf}
\usepackage{nccmath}
\usepackage{float}
\usepackage{subcaption}
\usepackage[noend]{algpseudocode}
\usepackage{algorithm2e}
\usepackage{multirow}
\usepackage{tabularx}
\usepackage{cite}
\usepackage{longtable}
\graphicspath{ {figures/} }
\usepackage{xargs}
\usepackage{orcidlink}

\allowdisplaybreaks

\IEEEoverridecommandlockouts

\begin{document}
\title{Numerical Performance of Different Formulations for Alternating Current Optimal Power Flow}

\author{Sayed Abdullah Sadat\textsuperscript{\Large \orcidlink{0000-0001-8290-6065}}~
        and~Kibaek Kim\textsuperscript{\Large \orcidlink{0000-0002-5820-6533}}
\thanks{This material is based upon work supported by the U.S. Department of Energy, Office of Science, under contract number DE-AC02-06CH11357.


Sayed Abdullah Sadat and Kibaek Kim are with the Mathematics and Computer Science Division, Argonne National Laboratory, The U.S. Department of Energy, Lemont, IL 60439, USA. (e-mail: sayed\_abdullah@ieee.org; kimk@anl.gov).

\noindent
U.S. Government work not protected by U.S. copyright.}}

\markboth{}%
{Shell \MakeLowercase{\textit{et al.}}: Bare Demo of IEEEtran.cls for IEEE Journals}

\maketitle
\thispagestyle{plain}
\pagestyle{plain}
\begin{abstract}
Alternating current optimal power flow (ACOPF) problems are nonconvex and nonlinear optimization problems.
Utilities and independent service operators (ISO) require  ACOPF  to be solved in almost real time. 
Interior point methods (IPMs) are one of the powerful methods for solving large-scale nonlinear optimization problems and are a suitable approach for solving ACOPF with large-scale real-world transmission networks. 
Moreover, the choice of the formulation is as important as choosing the algorithm for solving an ACOPF problem.  
In this paper, different ACOPF formulations with various linear solvers and the impact of employing box constraints are evaluated for computational viability and best performance when using IPMs. Different optimization structures are used in these formulations to model the ACOPF problem representing a range of  sparsity. The numerical experiments suggest that the least sparse ACOPF formulations with polar voltages yield the best computational results. Additionally, nodal injected models and current-based branch flow models are improved by enforcing box constraints. A wide range of test cases, ranging from 500-bus systems to 9591-bus systems, are used to verify the test results. 
\end{abstract}


%
\IEEEpeerreviewmaketitle

\section{Introduction}

\IEEEPARstart{P}{ower} system operation aims to deliver power to customers in a reliable and cost-effective manner. Optimal and reliable power dispatch is an optimization problem, which in its original form is referred to as the alternating current optimal power flow (ACOPF) problem \cite{ACOPFhistory,security}. ACOPF  is a nonlinear nonconvex optimization problem that is NP-hard.


Because of the computational difficulty of ACOPF, the majority of  grid operators use some form of a direct current optimal power flow (DCOPF) formulation \cite{stott2009dc}. The DCOPF linearization used in electricity market problems is based on a few assumptions that help convert the nonlinear ACOPF problem to a linear formulation \cite{sayed,sadat2021cascading}. Although the implementations enable solving optimal power flow within  practical time constraints, rough approximations of voltage and reactive power flow are applied, requiring the application of conservative transmission limits. This can most likely result in an inefficient and suboptimal operation of the grid. 

Over the past few decades, many methods have been proposed to solve the ACOPF problems. However, the majority of them overlook the impact of formulation on the performance of ACOPF solutions, especially when solving large-scale problems. Many methods have been proposed and tested for solving ACOPF, including a variety of convex relaxation techniques (e.g.,~\cite{lavaei2012zero,madani2015convex,low2014convex,hijazi2017convex,zhang2021tightness,sadat,SLP_sampath,sadat3}). 

In \cite{survey}, different optimal power flow formulations are presented and in \cite{park} the performance of three  ACOPF formulations is numerically analyzed. However, the analysis is limited to branch flow models, which  do not encompass the least sparse formulations for an ACOPF problem.


The past work presents a great deal of progress in solving ACOPF problems by offering different choices of the optimization algorithms. However, the performance of these algorithms can significantly be enhanced by the appropriate choice of the ACOPF formulation. None of the past work extensively explores the effect of formulation structure on the performance of their proposed algorithms, the impact of which could be significant, especially for large-scale problems.

This paper  evaluates the performance of different ACOPF formulations and identifies the best scalable ACOPF formulations. The study includes both branch flow and nodal injection models for an ACOPF problem. These formulations, while presenting the same problem, have different optimization structures. Generally, nodal injection formulations are least sparse in structure while branch flow formulation can be relatively sparse. The performance of these formulations is evaluated by using one of the most commonly used solvers, IPOPT, that implements the interior point method (IPM) with line search, which is typically used as a benchmark by most of the recent work on ACOPF algorithms \cite{kardos2020complete,castillo2013computational,lu2018optimization}.

In addition to sparsity, box constraints can  impact the performance of nonlinear programming. Although they are necessary for ensuring convergence in first-order methods for nonlinear programming \cite{bazara}, they can also impact the total number of iterations and computation time in interior point methods. 

The contributions of this paper can be summarized as follows:

\begin{enumerate}
\item Studying the impact of box constraints on ACOPF formulations using interior point method 
\item Evaluating the computational performance of different ACOPF formulations with respect to solution scalability, different voltage models, and sparsity structures
\item Identifying the best scalable ACOPF formulations and the factors contributing to better performance
\end{enumerate}

Similarly, different algorithms can have different performance behavior corresponding to different ACOPF formulations. The convergence pattern might differ when used in other methods. Thus, as a future work it is worth investigating the performance of different ACOPF formulations using different algorithms such as sequential linear programming (SLP), sequential quadratic programming (SQP), and alternating direction method of multipliers (ADMM); see, for example,~\cite{kim2021leveraging}.

The rest of the paper is organized as follows. Section II introduces optimal power flow, its different formulations, and box constraints. Section III presents the numerical experiments for evaluation of different ACOPF formulations and linear solvers. Section IV summarizes the conclusions of the study.

\section{Optimal Power Flow Formulations}
The ACOPF formulations developed in this work is based on the conventions used in \cite{sayed_unpublished} with the objective function \eqref{mainACOPF_obj}.
The formulation of interest in this paper is  as follows:
\begin{subequations}
\label{mainACOPF_obj}
\begin{align}
\underset{}{\text{min}}
& \sum_{\substack{g \in \mathcal{G}}} c_{2g} \cdot (\Re[s_g])^2 + c_{1g} \cdot \Re[s_g] + c_{0g}
\end{align}
\end{subequations}
Let $c_{2g},c_{1g},$ and $c_{0g}$ be the coefficients of quadratic cost function of generator $g$; let $s_g := p_g + j q_g$ be the apparent power output of generator $g$; let $v_n := v^r_n+jv^i_n = V_n e^{j\theta_n}$ be the voltage at bus $n$; let $S_{ij} := P_{ij}+jQ_{ij}$ be the line apparent power flow from node $i$ to $j$; let  $s_{d} := p_{d}+jq_{d}$ be the complex load; let $V_n^{Min}$ and $V_n^{Max}$ respectively be the lower and upper voltage bounds at bus $n$; let $I_{ij}^{Max}$ be the current thermal limit of line $\{i,j\}$; let $s_g^{Min} := p_g^{Min}+jq_g^{Min}, s_g^{Max}:=p_g^{Max}+jq_g^{Max}$ respectively be the complex generator $g$ lower and upper bounds; and let  $Y^L_{n}:=G^L_n-jB^L_n$ be the complex shunt element connected to bus $n$. Let $\mathcal{N}$ be the set of buses, $\mathcal{G}$ be the set of generators, and $\mathcal{K}$ be the set of branches in the network.

Power generation output can be separated into real and reactive power components, and therefore can be written as follows:

\begin{subequations}
\label{eq:ACOPF_pg_qg}
\begin{align}
& & &p_g^{Min} \leq p_g \leq p_g^{Max}, \; \forall g \in \mathcal{G}. \\
& & & q_g^{Min} \leq q_g \leq q_g^{Max}, \; \forall g \in \mathcal{G}.
\end{align}
\end{subequations}

\begin{figure}[!htb]
\centering
\includegraphics[scale=0.3]{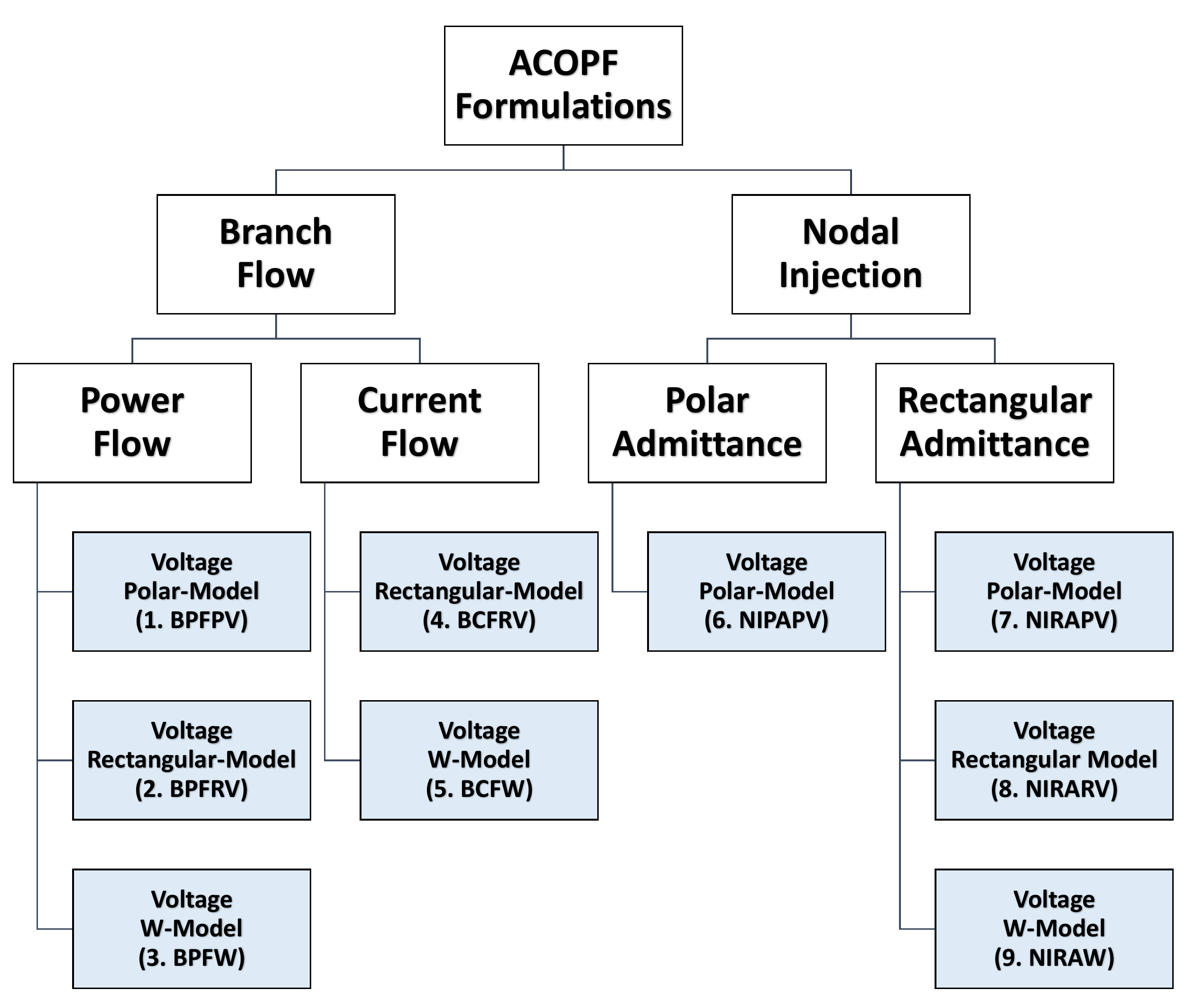}
\caption{Classification of different ACOPF formulations.}
\label{fig:ACOPFchart}
\end{figure}

An ACOPF problem in its complex form can be expanded into different solvable optimization formulations as shown in Fig. \ref{fig:ACOPFchart}. While these formulations represent the same original ACOPF problem, they have different optimization structures that can have different convergence patterns. The ACOPF formulations can be broadly classified as either branch flow models or nodal injection models.

In addition, the flow can be  represented by either power equations or current equations. Voltages are  represented by polar form, rectangular form, or $W$-model in modeling different ACOPF formulations.

\subsection{Branch Flow Model}

In a branch flow model ACOPF problem, a power or current flow in each branch is explicitly represented as variables. Using the power equations, one can model the branches as follows.

\begin{ceqn}
 \begin{align}
 \label{eq:Sij}
s_{ij}=v_i \cdot i_{ij}^*=v_i \cdot (y^{net}_{ij} \cdot v_i+Y_{ij} \cdot v_j)^*, \; \forall \{i,j\} \in 2\mathcal{K}
\end{align}
\end{ceqn}

Current flow equations can be written as 

\begin{ceqn}
 \begin{align}
 \label{eq:Iij}
i_{ij}=y^{net}_{ij} \cdot v_i+Y_{ij} \cdot v_j, \; \forall \{i,j\} \in 2\mathcal{K}.
\end{align}
\end{ceqn}

These branch variables (i.e., either $s_{ij}$ or $i_{ij}$), are then used in the node balance equations.

\subsection{Nodal Injection Model}
In this model, flows are not represented by explicit variables and instead are integrated in node power balance equations using admittance matrix and voltages as shown in \eqref{eq:ACOPF_NIM_si}. This approach reduces the sparsity of the ACOPF problem.

\begin{ceqn}
 \begin{align}
 \label{eq:ACOPF_NIM_si}
\sum_{\substack{g \in \mathcal{G}_n}}s_g - \sum_{\substack{d \in \mathcal{D}_n}}s_d = & v_n \cdot \sum_{\substack{k \in \mathcal{N}}}Y^*_{nk} \cdot v^*_k\nonumber \\
&+Y^L_{n}\cdot |v_n|^2, \; \forall n \in \mathcal{N}
\end{align}
\end{ceqn}

Three different approaches for modeling voltages in an ACOPF problem are proposed. 

\subsection{Voltage Polar Form}
In the voltage polar category, the voltages are modeled in polar form. Using the polar form of the voltages in ACOPF problem, one can model the branch power flow  as follows: 

\begin{subequations}
\label{eq:PFPV}
\begin{align}
& P_{ij} &&= g^{net}_{ij} \cdot V^2_i + (G_{ij} \cdot cos(\theta_i - \theta_j) \nonumber \\
& & & + B_{ij} \cdot sin(\theta_i - \theta_j) \cdot V_i \cdot V_j,\; \forall \{i,j\} \in 2\mathcal{K}. \\
& Q_{ij} &&= -b^{net}_{ij} \cdot V^2_i + (G_{ij} \cdot sin(\theta_i - \theta_j) \nonumber \\
& & & - B_{ij} \cdot cos(\theta_i - \theta_j) \cdot V_i \cdot V_j,\; \forall \{i,j\} \in 2\mathcal{K}.
\end{align}
\end{subequations}

Correspondingly, the node power balance can be written as  follows:

\begin{subequations}
\label{eq:PPV}
\begin{align}
& \sum_{\substack{g \in \mathcal{G}_n}}p_g - \sum_{\substack{d \in \mathcal{D}_n}}p_d = \sum_{\substack{k \in \mathcal{K}_n}}P_{nk}+G^L_{n}\cdot V_n^2, \; \forall n \in \mathcal{N}. \\
& \sum_{\substack{g \in \mathcal{G}_n}}q_g - \sum_{\substack{d \in \mathcal{D}_n}}q_d = \sum_{\substack{k \in \mathcal{K}_n}}Q_{nk}-B^L_{n}\cdot V_n^2, \; \forall n \in \mathcal{N}. 
\end{align}
\end{subequations}

For branch flow, the power thermal limit of transmission branches can be enforced as 

\begin{ceqn}
 \begin{align}
 \label{eq:BFPPV}
P^2_{ij} + Q^2_{ij} \leq (I^{Max}_{ij})^2 \cdot V^2_n, \; \forall \{i,j\} \in 2\mathcal{K}. 
\end{align}
\end{ceqn}

Similarly, for the nodal injection models with polar admittance, the nodal power balance can be written as follows:

\begin{subequations}
\label{eq:NIPPAPV}
\begin{align}
& \sum_{\substack{g \in \mathcal{G}_n}}p_g - \sum_{\substack{d \in \mathcal{D}_n}}p_d &&= \sum_{\substack{k \in \mathcal{N}}}|Y_{nk}| \cdot V_n \cdot V_k \cdot cos(\theta_n -\theta_k\nonumber \\
& & & -\theta^Y_{nk})+G^L_{n}\cdot V^2_n, \; \forall n \in \mathcal{N}.\\
& \sum_{\substack{g \in \mathcal{G}_n}}q_g - \sum_{\substack{d \in \mathcal{D}_n}}q_d &&= \sum_{\substack{k \in \mathcal{N}}}|Y_{nk}| \cdot V_n \cdot V_k \cdot sin(\theta_n-\theta_k\nonumber \\
& & & -\theta^Y_{nk})-B^L_{n}\cdot V^2_n, \; \forall n \in \mathcal{N}.
\end{align}
\end{subequations}

Similarly, for the nodal injection model with rectangular admittance, the nodal power balance can be written as follows:

\begin{subequations}
\label{eq:NIPRAPV}
\begin{align}
& \sum_{\substack{g \in \mathcal{G}_n}}p_g - \sum_{\substack{d \in \mathcal{D}_n}}p_d &&= \sum_{\substack{k \in \mathcal{N}}}V_n \cdot V_k \cdot (G_{nk} \cdot cos(\theta_n-\theta_k)\nonumber \\
& & & +B_{nk} \cdot sin(\theta_n-\theta_k))\nonumber \\
& & & +G^L_{n}\cdot V^2_n, \; \forall n \in \mathcal{N}.\\
& \sum_{\substack{g \in \mathcal{G}_n}}q_g - \sum_{\substack{d \in \mathcal{D}_n}}q_d &&= \sum_{\substack{k \in \mathcal{N}}}V_n \cdot V_k \cdot (G_{nk} \cdot sin(\theta_n-\theta_k)\nonumber \\
& & & -B_{nk} \cdot cos(\theta_n-\theta_k))\nonumber \\
& & & -B^L_{n}\cdot V^2_n, \; \forall n \in \mathcal{N}.
\end{align}
\end{subequations}

For nodal injection models, the current thermal limit of transmission branches can be enforced as follows:

\begin{subequations}
\label{eq:NIIPV}
\begin{align}
& 2 \cdot |y^{net}_{ij}| \cdot |Y_{ij}| \cdot V_i\cdot V_j \cdot [\nonumber \\
& cos(\theta_n -\theta_k+\theta^{yY}_{nk})] \nonumber \\
& +|y^{net}_{ij}|^2 \cdot V^2_i+|Y_{ij}|^2 \cdot V^2_j
\ \leq (I^{Max}_{ij})^2,\; \forall \{i,j\} \in 2\mathcal{K}.
\end{align}
\end{subequations}

The voltage magnitude limits can be enforced by using the following:
\begin{ceqn}
 \begin{align}
 \label{eq:V_limits}
V_n^{Min} \leq V_n \leq V_n^{Max}, \; \forall n \in \mathcal{N}. 
\end{align}
\end{ceqn}

\subsection{Voltage Rectangular Form}
In the voltage rectangular category, voltages are modeled in rectangular form. Using the rectangular form of the voltage in the ACOPF problem, the branch power flow can be modeled as follows:

\begin{subequations}
\label{eq:PFRV}
\begin{align}
& P_{ij}&&=g^{net}_{ij} \cdot [(v^r_i)^2+(v^i_i)^2] + G_{ij} \cdot v_i^r \cdot v_j^r + B_{ij} \cdot v_i^i \cdot v_j^r \nonumber \\
& & & + G_{ij} \cdot v_i^i \cdot v_j^i - B_{ij} \cdot v_i^r \cdot v_j^i, \; \forall \{i,j\} \in 2\mathcal{K}. \\
& Q_{ij}&&=-b^{net}_{ij} \cdot [(v^r_i)^2+(v^i_i)^2] + G_{ij} \cdot v_i^i \cdot v_j^r - B_{ij} \cdot v_i^r \cdot v_j^r \nonumber \\
& & & - G_{ij} \cdot v_i^r \cdot v_j^i - B_{ij} \cdot v_i^i \cdot v_j^i, \; \forall \{i,j\} \in 2\mathcal{K}. 
\end{align}
\end{subequations}

Correspondingly, the node power balance can be written as follows:

\begin{subequations}
\label{eq:PRV}
\begin{align}
& \sum_{\substack{g \in \mathcal{G}_n}}p_g - \sum_{\substack{d \in \mathcal{D}_n}}p_d=\sum_{\substack{k \in \mathcal{K}_n}}P_{nk}+G^L_{n}\cdot v^2_n, \; \forall n \in \mathcal{N}. \\
& \sum_{\substack{g \in \mathcal{G}_n}}q_g - \sum_{\substack{d \in \mathcal{D}_n}}q_d=\sum_{\substack{k \in \mathcal{K}_n}}Q_{nk}-B^L_{n}\cdot v^2_n, \; \forall n \in \mathcal{N}. 
\end{align}
\end{subequations}
where $v_n^2 = (v^r_n)^2+(v^i_n)^2$.

For branch power flow, the thermal limit of transmission branches can be enforced as follows:

\begin{ceqn}
 \begin{align}
 \label{eq:BFPRV}
P^2_{ij} + Q^2_{ij} \leq (I^{Max}_{ij})^2 \cdot [(v^r_n)^2+(v^i_n)^2], \; \forall \{i,j\} \in 2\mathcal{K}.
\end{align}
\end{ceqn}

Using the rectangular form of the voltage in the ACOPF problem, the branch current flow can be modeled as follows:

\begin{subequations}
\label{eq:CFRV}
\begin{align}
& i^r_{ij}&&=g^{net}_{ij} \cdot v^r_i-b^{net}_{ij} \cdot v^i_i+G_{ij} \cdot v^r_j\nonumber \\
& & &-B_{ij} \cdot v^i_j, \; \forall \{i,j\} \in 2\mathcal{K}.   \\
& i^i_{ij}&&=g^{net}_{ij} \cdot v^i_i+b^{net}_{ij} \cdot v^r_i+G_{ij} \cdot v^i_j\nonumber \\
& & &+B_{ij} \cdot v^r_j, \; \forall \{i,j\} \in 2\mathcal{K}. 
\end{align}
\end{subequations}

Correspondingly, the node power balance can be written as follows:

\begin{subequations}
\label{eq:CRV}
\begin{align}
& \sum_{\substack{g \in \mathcal{G}_n}}p_g - \sum_{\substack{d \in \mathcal{D}_n}}p_d &&=\sum_{\substack{\{n,k\} \in \mathcal{K}_n^f \bigcup \mathcal{K}_n^t}}(v^r_n \cdot i^r_{nk}+v^i_n \cdot i^i_{nk})\nonumber \\
& & &+G^L_{n}\cdot [(v^r_n)^2+(v^i_n)^2], \; \forall n \in \mathcal{N}. \\
& \sum_{\substack{g \in \mathcal{G}_n}}q_g - \sum_{\substack{d \in \mathcal{D}_n}}q_d && =\sum_{\substack{\{n,k\} \in \mathcal{K}_n^f \bigcup \mathcal{K}_n^t}}(v^i_n \cdot i^r_{nk}+v^r_n \cdot i^i_{nk})\nonumber \\
& & &-B^L_{n}\cdot [(v^r_n)^2+(v^i_n)^2], \; \forall n \in \mathcal{N}, 
\end{align}
\end{subequations}
where $i_{ij} = i^r_{ij}+ji^i_{ij}$ is the complex line current flow from node $i$ to $j$.

For branch current flow, the thermal limit of transmission branches can be enforced as follows:

\begin{ceqn}
 \begin{align}
 \label{eq:BFCRV}
(i^r_{ij})^2 + (i^i_{ij})^2 \leq (I^{Max}_{ij})^2 , \; \forall \{i,j\} \in 2\mathcal{K}.
\end{align}
\end{ceqn}

Similarly, for the nodal injection model with rectangular admittance, the nodal power balance can be written as \eqref{eq:NIPRARV}.

\begin{subequations}
\label{eq:NIPRARV}
\begin{align}
& \sum_{\substack{g \in \mathcal{G}_n}}p_g - \sum_{\substack{d \in \mathcal{D}_n}}p_d &&= \sum_{\substack{k \in \mathcal{N}}}
(G_{nk} \cdot v^r_n \cdot v^r_k + B_{nk} \cdot v^i_n \cdot v^r_k \nonumber \\
& & & + G_{nk} \cdot v^i_n \cdot v^i_k - B_{nk} \cdot v^r_n \cdot v^i_k) \nonumber \\
& & &+G^L_{n}\cdot [(v^r_n)^2+(v^i_n)^2], \; \forall n \in \mathcal{N}.\\
& \sum_{\substack{g \in \mathcal{G}_n}}q_g - \sum_{\substack{d \in \mathcal{D}_n}}q_d &&= \sum_{\substack{k \in \mathcal{N}}}
(G_{nk} \cdot v^i_n \cdot v^r_k - B_{nk} \cdot v^r_n \cdot v^r_k \nonumber \\
& & &- G_{nk} \cdot v^r_n \cdot v^i_k - B_{nk} \cdot v^i_n \cdot v^i_k) \nonumber \\
& & &-B^L_{n}\cdot [(v^r_n)^2+(v^i_n)^2], \; \forall n \in \mathcal{N}.
\end{align}
\end{subequations}

For nodal injection models, the thermal limit of transmission branches can be enforced as follows:

\begin{subequations}
\label{eq:NIIRV}
\begin{align}
&(g^{net}_{ij} \cdot v^r_i-b^{net}_{ij} \cdot v^i_i \nonumber \\
& +G_{ij} \cdot v^r_j-B_{ij} \cdot v^i_j)^2 \nonumber \\
& + (g^{net}_{ij} \cdot v^i_i+b^{net}_{ij} \cdot v^r_i \nonumber \\
& +G_{ij} \cdot v^i_j+B_{ij} \cdot v^r_j)^2 
\ \leq (I^{Max}_{ij})^2,\; \forall \{i,j\} \in 2\mathcal{K}.
\end{align}
\end{subequations}

The voltage magnitude limits can be enforced by using 
\begin{ceqn}
 \begin{align}
 \label{eq:v_limits}
(V_n^{Min})^2 \leq (v^r_n)^2+(v^i_n)^2 \leq (V_n^{Max})^2, \; \forall n \in \mathcal{N}.
\end{align}
\end{ceqn}

\subsection{Voltage W-Model} 
In the voltage W-Model category, we replace the nonlinearities of the terms containing voltages with W-Model, thereby making the problem more sparse. The $W$-matrix and the bounds on its elements are defined in \eqref{eq:W_limits}. The branch power flow can be modeled as follows:

\begin{subequations}
\label{eq:PFW}
\begin{align}
& P_{ij}=g^{net}_{ij} \cdot W^d_i + G_{ij} \cdot W^r_{ij} - B_{ij} \cdot W^i_{ij}, \; \forall \{i,j\} \in 2\mathcal{K}. \\
& Q_{ij}=-b^{net}_{ij} \cdot W^d_i - G_{ij} \cdot W^i_{ij} - B_{ij} \cdot W^r_{ij}, \; \forall \{i,j\} \in 2\mathcal{K}.
\end{align}
\end{subequations}

Correspondingly, the node power balance can be written as follows:

\begin{subequations}
\label{eq:PW}
\begin{align}
& \sum_{\substack{g \in \mathcal{G}_n}}p_g - \sum_{\substack{d \in \mathcal{D}_n}}p_d =\sum_{\substack{k \in \mathcal{K}_n}}P_{nk}+G^L_{n}\cdot W^d_n, \; \forall n \in \mathcal{N} \\
& \sum_{\substack{g \in \mathcal{G}_n}}q_g - \sum_{\substack{d \in \mathcal{D}_n}}q_d =\sum_{\substack{k \in \mathcal{K}_n}}Q_{nk}-B^L_{n}\cdot W^d_n, \; \forall n \in \mathcal{N}. 
\end{align}
\end{subequations}

For branch flow, the power thermal limit of transmission branches can be enforced as follows:

\begin{ceqn}
 \begin{align}
 \label{eq:BFPW}
P^2_{ij} + Q^2_{ij} \leq (I^{Max}_{ij})^2 \cdot W^d_n, \; \forall \{i,j\} \in 2\mathcal{K}. 
\end{align}
\end{ceqn}

The node power balance for the branch current flow can be modeled as follows:

\begin{subequations}
\label{eq:CW}
\begin{align}
& \sum_{\substack{g \in \mathcal{G}_n}}p_g - \sum_{\substack{d \in \mathcal{D}_n}}p_d &&=\sum_{\substack{\{n,k\} \in \mathcal{K}_n^f \bigcup \mathcal{K}_n^t}}(v^r_n \cdot i^r_{nk}+v^i_n \cdot i^i_{nk}) \nonumber \\
& & & +G^L_{n}\cdot W^d_n, \; \forall n \in \mathcal{N}. \\
& \sum_{\substack{g \in \mathcal{G}_n}}q_g - \sum_{\substack{d \in \mathcal{D}_n}}q_d &&=\sum_{\substack{\{n,k\} \in \mathcal{K}_n^f \bigcup \mathcal{K}_n^t}}(v^i_n \cdot i^r_{nk}+v^r_n \cdot i^i_{nk}) \nonumber \\
& & & -B^L_{n}\cdot W^d_n, \; \forall n \in \mathcal{N}.
\end{align}
\end{subequations}

Similarly, for the nodal injection models with rectangular admittance the nodal power balance can be written as follows:

\begin{subequations}
\label{eq:NIPRAW}
\begin{align}
& \sum_{\substack{g \in \mathcal{G}_n}}p_g - \sum_{\substack{d \in \mathcal{D}_n}}p_d &&= \sum_{\substack{k \in \mathcal{N}}}
(G_{nk} \cdot W^r_{nk} - B_{nk} \cdot W^i_{nk}) \nonumber \\
& & & +G^L_{n}\cdot W^d_n, \; \forall n \in \mathcal{N}.\\
& \sum_{\substack{g \in \mathcal{G}_n}}q_g - \sum_{\substack{d \in \mathcal{D}_n}}q_d &&= -\sum_{\substack{k \in \mathcal{N}}}
(G_{nk} \cdot W^i_{nk} + B_{nk} \cdot W^r_{nk}) \nonumber \\
& & & -B^L_{n}\cdot W^d_n, \; \forall n \in \mathcal{N}.
\end{align}
\end{subequations}

For nodal injection models, the current thermal limit of transmission branches can be enforced as follows:

\begin{subequations}
\label{eq:NIIW}
\begin{align}
& 2 \cdot (g^{net}_{ij} \cdot G_{ij} + b^{net}_{ij} \cdot B_{ij} ) \cdot W^r_{ij} \nonumber \\
& +2 \cdot (b^{net}_{ij} \cdot G_{ij} - g^{net}_{ij} \cdot B_{ij} ) \cdot W^i_{ij} \nonumber \\
& +|y^{net}_{ij}|^2 \cdot W^d_i + |Y_{ij}|^2 \cdot W^d_j \leq (I^{Max}_{ij})^2,\; \forall \{i,j\} \in 2\mathcal{K}.
\end{align}
\end{subequations}

The W-Model and the limits can be formed by using the voltage relationship  as follows:
\begin{subequations}
 \label{eq:W_limits}
 \begin{align}
& (V_n^{Min})^2 \leq W^d_{n} \leq (V_n^{Max})^2, \; \forall n \in \mathcal{N}. \label{eq:Wd_limits} \\
& (W^r_{ij})^2 + (W^i_{ij})^2 = W^d_{i} \cdot W^d_{j}, \; \forall \{i,j\} \in 2\mathcal{K}. \\
& W^d_{n} = (v^r_n)^2+(v^i_n)^2, \; \forall n \in \mathcal{N}.\label{eq:Wd_vrvi} \\
& W^r_{ij} = v^r_i \cdot v^r_j+v^i_i \cdot v^i_j, \; \forall n \in \mathcal{K}.\\
& W^i_{ij} = v^r_i \cdot v^i_j+v^i_i \cdot v^r_j, \; \forall n \in \mathcal{K}.
\end{align}
\end{subequations}

\begin{table}[h!]
\begin{center}
\caption{Analysis of the optimization structure of different ACOPF formulations with objective function \eqref{mainACOPF_obj} without considering the box constraints. Typically, for a transmission system $\mathcal{K} > \mathcal{N} > \mathcal{G}$.}
\label{table:formulations}
\begin{tabular}{|c|l|c|c|}
\hline \hline
  \multirow{2}{*}{Sno.} &\multicolumn{1}{c|}{\multirow{1}{*}{ACOPF}} &\multirow{2}{*}{Constraints} &Nonlinear \\ 
  &\multicolumn{1}{c|}{\multirow{1}{*}{Formulations}} & &Constraints  \\ \hline
  1&BPFPV &\eqref{eq:ACOPF_pg_qg},\eqref{eq:PFPV},\eqref{eq:PPV},\eqref{eq:BFPPV},\eqref{eq:V_limits} &$2\mathcal{N}+6\mathcal{K}$ \\ \hline
  2&BPFRV &\eqref{eq:ACOPF_pg_qg},\eqref{eq:PFRV},\eqref{eq:PRV},\eqref{eq:BFPRV},\eqref{eq:v_limits} &$4\mathcal{N}+6\mathcal{K}$ \\ \hline 
  3&BPFW  &\eqref{eq:ACOPF_pg_qg},\eqref{eq:PFW},\eqref{eq:PW},\eqref{eq:BFPW},\eqref{eq:W_limits}  &$\mathcal{N}+6\mathcal{K}$ \\ \hline
  4&BCFRV &\eqref{eq:ACOPF_pg_qg},\eqref{eq:CFRV},\eqref{eq:CRV},\eqref{eq:BFCRV},\eqref{eq:v_limits}   &$4\mathcal{N}+2\mathcal{K}$ \\ \hline
  5&BCFW &\eqref{eq:ACOPF_pg_qg},\eqref{eq:CFRV},\eqref{eq:CW},\eqref{eq:BFCRV},\eqref{eq:Wd_limits},\eqref{eq:Wd_vrvi}  &$3\mathcal{N}+6\mathcal{K}$ \\ \hline
  6&NIPAPV &\eqref{eq:ACOPF_pg_qg},\eqref{eq:NIPPAPV},\eqref{eq:NIIPV},\eqref{eq:V_limits}  &$2\mathcal{N}+2\mathcal{K}$ \\ \hline
  7&NIRAPV &\eqref{eq:ACOPF_pg_qg},\eqref{eq:NIPRAPV},\eqref{eq:NIIPV},\eqref{eq:V_limits}  &$2\mathcal{N}+2\mathcal{K}$ \\ \hline
  8&NIRARV &\eqref{eq:ACOPF_pg_qg},\eqref{eq:NIPRARV},\eqref{eq:NIIRV},\eqref{eq:v_limits}  &$4\mathcal{N}+2\mathcal{K}$ \\ \hline
  9&NIRAW &\eqref{eq:ACOPF_pg_qg},\eqref{eq:NIPRAW},\eqref{eq:NIIW},\eqref{eq:W_limits}   &$\mathcal{N}+4\mathcal{K}$ \\ \hline
  \hline
\end{tabular}
\end{center}
\end{table}

\subsection{Box Constraints}
\label{S:box}

\begin{figure}[h!]
\centering
\includegraphics[scale=0.33]{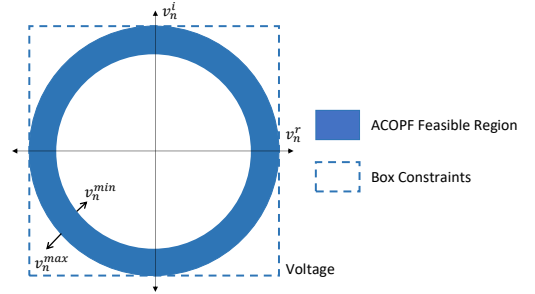}
\caption{Illustration of box constraints for the nodal voltage at bus $n$.}
\label{fig:box}
\end{figure}

Additional box constraints can be added to some of the ACOPF formulations.
While redundant to the formulations, the additional constraints can often be useful to better define the feasible solution region of optimization problem.
In an ACOPF problem, the box constraints can be applied to complex parameters based on the magnitude of the parameter (e.g., nodal voltages, branch power flows, and branch current flows). We study the impact of the additional box constraints with respect to the computational performance of these formulations. An illustrative example of box constraints is shown in Fig. \ref{fig:box} for a nodal voltage at bus $n$.

\section{Numerical Experiments}

All numerical experiments were performed on an 8-core Intel(R) i9-9980HK CPU @ 2.40 GHz with 32 GB RAM system. The models were implemented in Julia v1.5.3 using Ipopt.jl v0.6.5, and JuMP.jl v0.21.5 \cite{jump} in a Linux Ubuntu operating system. 
In the IPOPT solver, we experiment with different linear solvers: MUMPS, MA27, MA57, MA86, MA97, and PARDISO \cite{pardiso}. Moreover, each linear solver is limited to using one core. 

\begin{figure}[h!]
\centering
\includegraphics[scale=0.49]{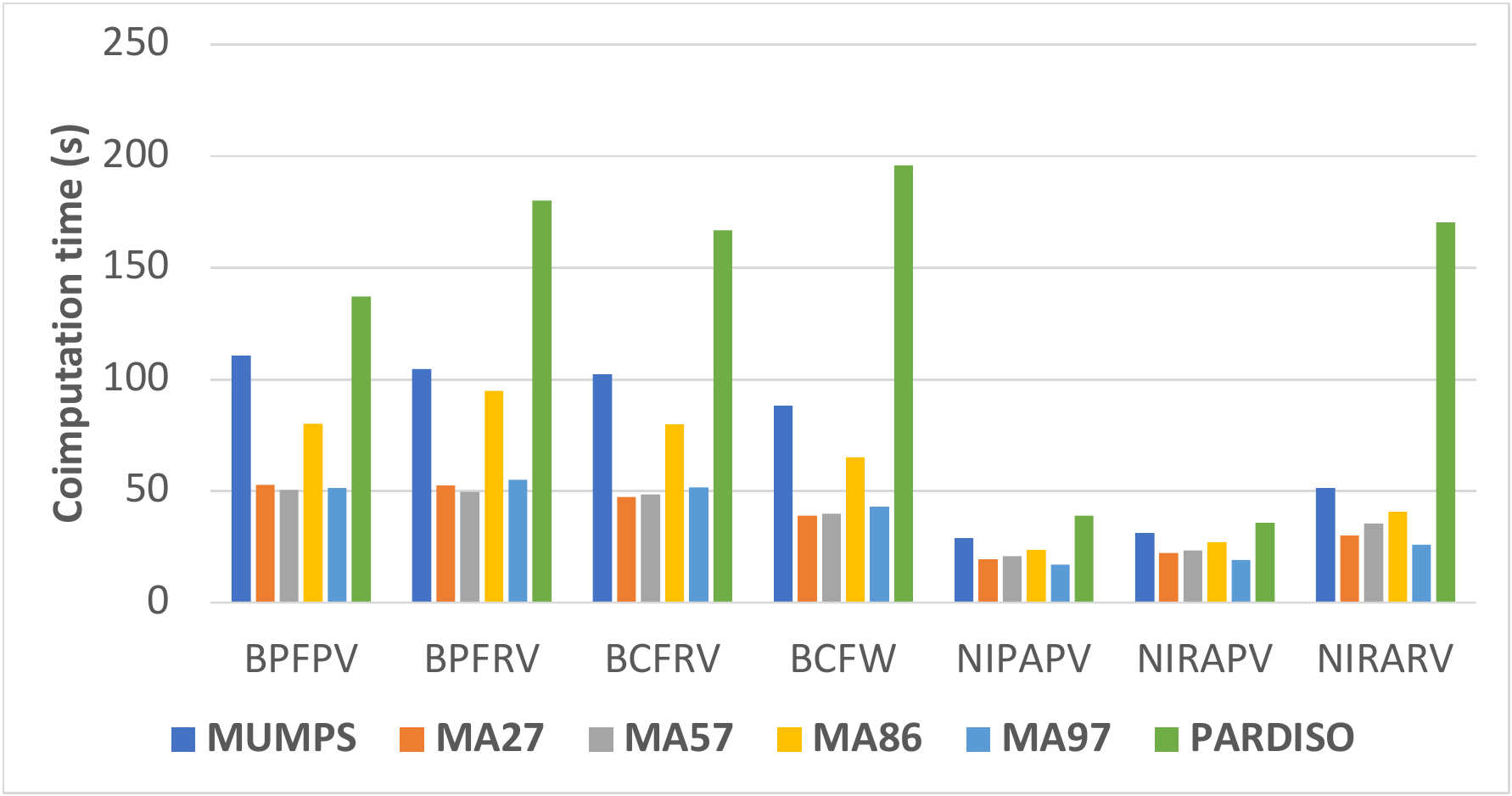}
\caption{Computation time of different ACOPF formulations using different linear solvers for solving Network\_84R-100 using the IPOPT solver.}
\label{fig:formulations_solver}
\end{figure}

The performance of different ACOPF formulations using different linear solvers is shown in Fig. \ref{fig:formulations_solver}. The IPOPT solver was not able to solve the full $W$-models (i.e., BPFW and NIRAW) giving the error of ``Too few degrees of freedom." The results show that smaller, dense formulations (NIPAPV, NIRAPV)  yield better performance compared with other formulations. The results also show that among the formulations the nodal injection polar voltage models (i.e., NIPAPV and NIRAPV) converge faster. However, among the other sparse formulations the current-based flow models, particularly diagonal $W$-model (BCFW), report the smallest solution time. 

\begin{figure}[!htb]
\centering
\includegraphics[scale=0.285]{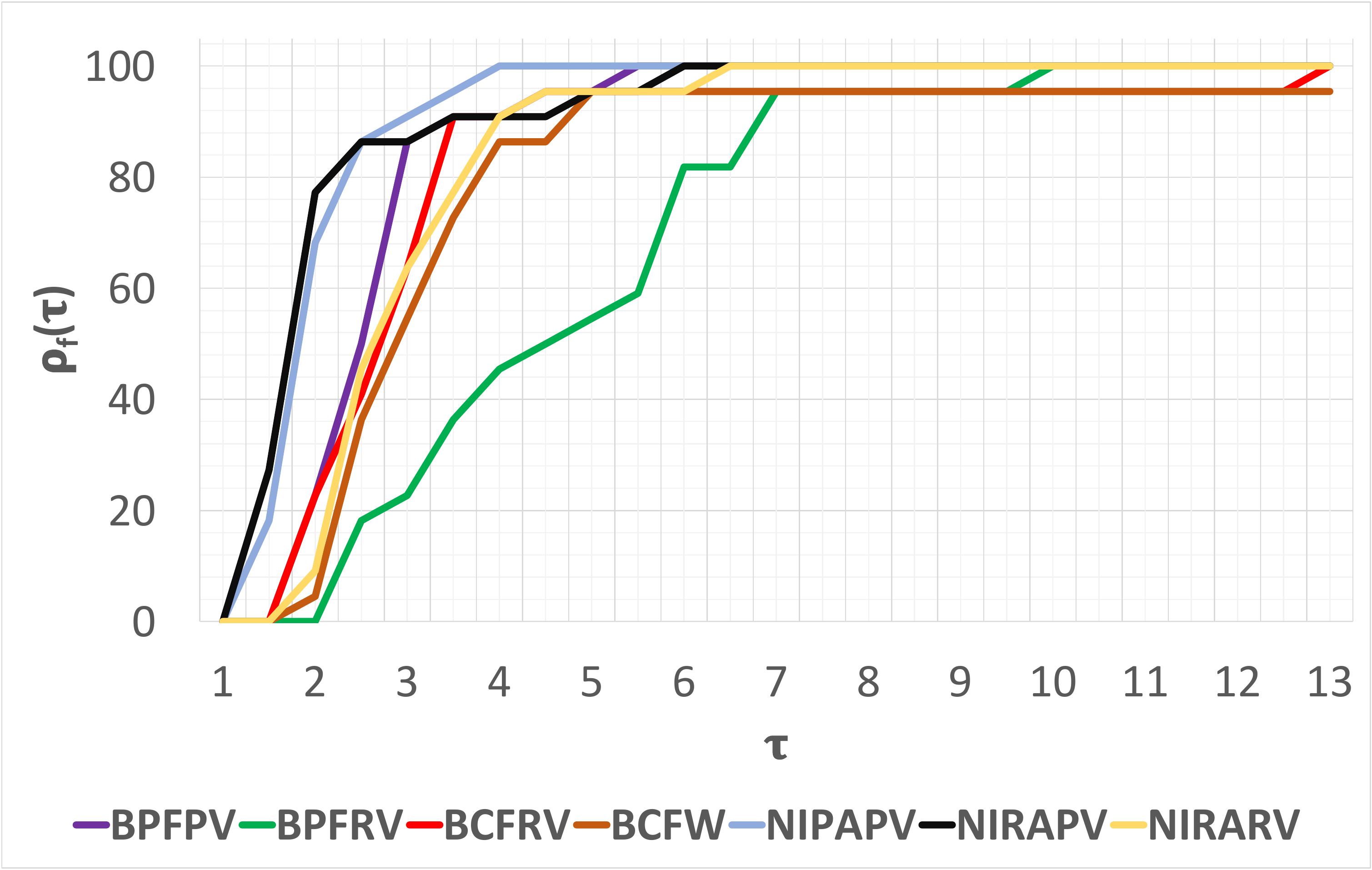}
\caption{Performance profile of different ACOPF formulations using IPOPT solver for different scenarios of a wide range of test cases shown in Table \ref{table:test_cases} \cite{sayed_unpublished}.}
\label{fig:performance}
\end{figure}

Figure \ref{fig:performance} shows the performance profile of the ACOPF formulations plotted as described in \cite{sayed_unpublished} by solving a wide range of different scenarios of test cases shown in Table \ref{table:test_cases}. 

\begin{table}[h!]
\begin{center}
\caption{Size of different test networks used in the numerical experiments \cite{GOch1, sayed_unpublished}.}
\label{table:test_cases}
\begin{tabular}{|c|c|c|c|}
\hline \hline
   \multirow{1}{*}{Networks} &\multirow{1}{*}{Buses} &\multirow{1}{*}{Gens} &\multirow{1}{*}{Branches} \\ \hline \hline
Network\_01R-040 &500 &90 &594  \\ \hline
Network\_03R-050 &793 &210 &912  \\ \hline
Network\_05R-018 &2000 &384 &3639 \\ \hline
Network\_07R-122 &2312 &617 &3013 \\ \hline
Network\_81R-050 &3288 &379 &4871 \\ \hline
Network\_84R-100 &4601 &410 &7304 \\ \hline
Network\_09R-070 &4918 &1340 &6727 \\ \hline
Network\_12R-030 &9591 &365 &15915  \\ \hline

\hline  \hline
\end{tabular}
\end{center}
\end{table}

\begin{figure*}[!htb]
\centering
\includegraphics[scale=0.49]{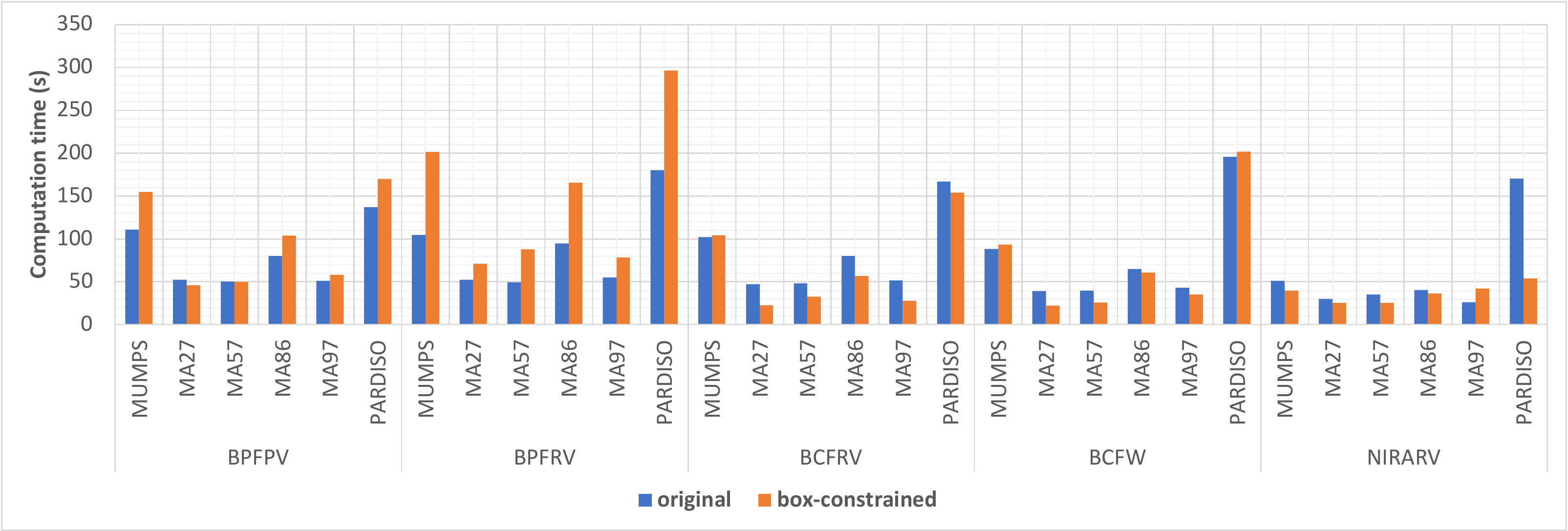}
\caption{Comparison of the computation time for solving different ACOPF formulations with and without box constraints of Network\_84R-100 by using the IPOPT solver.}
\label{fig:solver}
\end{figure*}

\begin{figure}[!htb]
\centering
\includegraphics[scale=0.51]{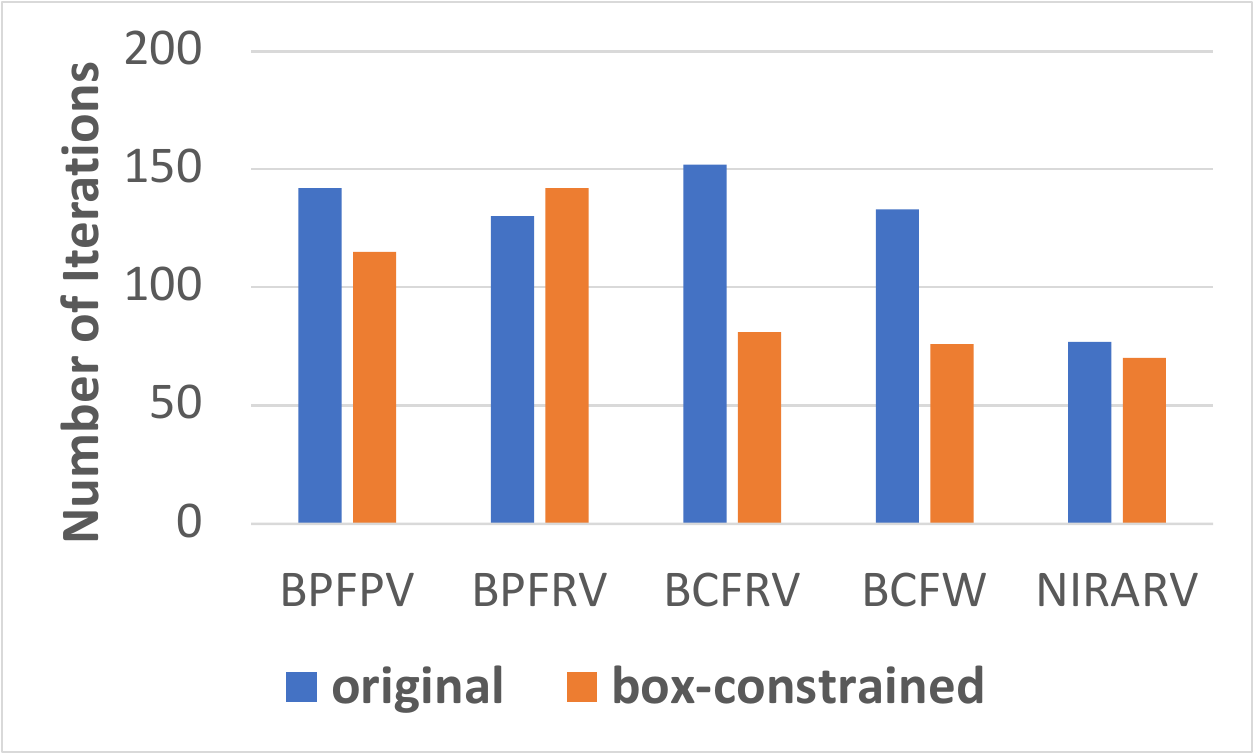}
\caption{Comparison of the iterations number for solving different ACOPF formulations with and without box constraints of Network\_84R-100 by using IPOPT solver.}
\label{fig:iterations}
\end{figure}

We apply the box constraints (see Section~\ref{S:box}) to the formulations that contain different optimization variables without bounds. These formulations include BPFPV, BPFRV, BCFRV, BCFW, and NIRARV. In Fig. \ref{fig:solver}  the computation time of the original ACOPF formulations (without box constraints)  with different linear solvers is compared with those with the box constraints. The results show that the box constraints affect the performance of ACOPF formulations and improve the performance of the BCFRV, BCFW, and NIRARV formulations. As shown in Fig. \ref{fig:iterations}, the number of iterations for BCFRV and BCFW were halved when box constraints are added, which may lead to the reduction in solution time for certain linear solvers.

\section{Conclusion}
This paper presents different ACOPF formulations, each representing a unique optimization structure and sparsity. The work includes an explicit investigation of the impact of these ACOPF formulations on the performance of ACOPF algorithms. The numerical experiment in this paper considers IPM implemented to evaluate these formulations by using the IPOPT solver. The results show that the choice of formulation, which is often overlooked, can significantly impact the computation time of an algorithm, especially for large-scale networks. Numerical results in the paper suggest a consistently superior performance by the Nodal Injection Polar Admittance Polar Voltage (NIPAPV) and Nodal Injection Rectangular Admittance Polar Voltage (NIRAPV) formulations among the different ACOPF formulations studied in this paper across a wide range of test networks. The complexity of the subproblem of these two formulations compared with the other formulations has also consistently remained the least. Box constraints can impact the performance of ACOPF formulations and improve the performance of BCFRV, BCFW, and NIRARV. It is worth investigating the performance of different ACOPF formulations with other algorithms such as SLP, SQP, and ADMM.

\ifCLASSOPTIONcaptionsoff
  \newpage
\fi

\bibliographystyle{IEEEtran}
\bibliography{IEEEabrv,anl_con_1}

\noindent\fbox{\parbox{0.47\textwidth}{
The submitted manuscript has been created by UChicago Argonne, LLC, Operator of Argonne National Laboratory (``Argonne''). Argonne, a U.S. Department of Energy Office of Science laboratory, is operated under Contract No. DE-AC02-06CH11357. The U.S. Government retains for itself, and others acting on its behalf, a paid-up nonexclusive, irrevocable worldwide license in said article to reproduce, prepare derivative works, distribute copies to the public, and perform publicly and display publicly, by or on behalf of the Government. The Department of Energy will provide public access to these results of federally sponsored research in accordance with the DOE Public Access Plan (http://energy.gov/downloads/doe-public-access-plan).}
}

\end{document}